\documentclass[graybox]{svmult}

\usepackage{type1cm}        
\usepackage{makeidx}         
\usepackage{graphicx}        
\usepackage{multicol}        
\usepackage[bottom]{footmisc}

\usepackage{newtxtext}   %
\usepackage[varvw]{newtxmath}       

\makeindex             

\def\vu{ {{\bf u}} }
\def\vE{ {{\bf E}} }
\def\vT{ {{\bf T}} }
\def\vx{ {{\bf x}} }
\def\v0{ {{\bf 0}} }

\begin{document}

\title*{A high order block finite difference, error inhibiting scheme for the transport equation}
\author{Adi Ditkowski\and Anne Le Blanc\and Chi-Wang Shu}
\institute{Adi Ditkowski \at Tel Aviv University, Tel Aviv 69978,Israel \email{adid@tauex.tau.ac.il}, \and Anne Le Blanc \at Tel Aviv University, Tel Aviv 69978,Israel, \email{anneleb@tauex.tau.ac.il}\and Chi-Wang Shu \at Brown University, Providence, RI 02912, USA, \email{chi-wang\_ shu@brown.edu}
	}

%
%
\maketitle

\abstract{We propose a block finite difference, error inhibiting scheme that is 
fourth-order accurate for short to moderate times and has a six-order convergence rate for long times. This scheme outperforms the standard fourth-order Finite Difference scheme. We also demonstrate that the proposed scheme is a particular type of nodal-based Discontinuous Galerkin method with $p=1$.}

\section{Introduction}\label{sec1}

The main theme of this manuscript is the construction and analysis of error-inhibiting block finite difference schemes for initial boundary value problems for PDEs. Applications to hyperbolic problems, namely the pure transport problem, are emphasized in this text.\\

Long-time numerical integration has become a significant issue in computational electrodynamics and fluid dynamics. A scheme might be stable in the classical sense (Lax stability), while its error grows exponentially in time \cite{carpenter}. Typically, the primary source of the error is the linearly growing phase error. The source of the phase error is that for each mode (frequency), the velocity of the numerical wave is slightly different from the analytical one. Therefore, reducing this phase error as much as possible is critical.\\

The standard procedure to build Finite Difference (FD) schemes is first, constructing a \textit{consistent} approximation for the given problem,  and secondly, its \textit{stability}, which is the discrete analog of well-posedness \cite{gustafsson1995time}, is proven. Lax's equivalence theorem \cite{richtmyer} ensures the convergence of the scheme, meaning that for a finite final time $T$, the numerical solution converges to the analytic solution of the problem as the mesh size goes to zero \cite{morton}. In the case of semi-discrete schemes, which is the class of schemes included in our work, stability is equivalent to the semi-boundedness of right hand side of the semi-discrete scheme. It implies that all eigenvalues of the coefficient matrix of the Ordinary Differential Equation (ODE) system have real parts bounded from above by a constant \cite{gustafsson1995time}, as well as a full set of eigenvectors that remain independent as the mesh is refined. However, this analysis is not refined enough for our purpose.\\

We focus on the different aspects involved in constructing and analyzing an efficient numerical scheme called Block-Finite-Difference, derived from the classical FD method.
The BFD schemes derived for the Heat equation (see \cite{ditkowski2020error}, \cite{leblanc2020error} for further details) relied on the inherent dissipation of the diffusion operator. This dissipation, as well as the post-processing procedure, caused the damping of the high-frequency elements of the error. The transport problem has no dissipation. Therefore, we rely on other mechanisms to manipulate the dynamics in which the error is accumulated in time.\\

The main idea behind the derivation of these BFD schemes for the transport problem is to manipulate the interaction between the truncation error and the dynamics of the scheme such that it inhibits the error from accumulating. In this manuscript, we derived a scheme with a third-order truncation error. However, the leading term of the error is of third-order, bounded in time, which can be eliminated in a post-processing stage to get fourth-order, bounded in time term, and the six-order phase error. Thus, this is a third-order or fourth-order after post-processing scheme for short to moderate times, while it is a six-order method for long times.\\

This manuscript is constructed as follows: in Section \ref{sec2}
 we construct the block finite difference scheme, analyze its stability and accuracy, and present several numerical examples to demonstrate its properties. Section \ref{sec:DG} shows that this scheme can be derived as a nodal-based DG method by adding extra penalty terms. Conclusions are given in Section \ref{sec:Conclusions}

\section{Block Finite Difference Methods for the Transport Problem}\label{sec2}

Consider the following problem
\begin{equation}
\label{eq1}
\left\lbrace 
\begin{array}{ll}

\dfrac{\partial u}{\partial t}+ \dfrac{\partial u}{\partial x } =0\mbox{ , }x\in \left( 0,2\pi \right) \mbox{ , }t\geq 0\\
u(x,0)=f(x)\\

\end{array}
  \right.
\end{equation}
with periodic boundary conditions.


\noindent We solve this problem over the following grid:
\begin{equation}
\; x_{j-1/4}=x_{j}-h/4\;,\;x_{j+1/4}=x_{j}+h/4\;,\; h=2\pi/N , \; j=1,...,N
\end{equation}
where 
\begin{equation}
x_{j}=h(j-1)+\frac{h}{2}
\end{equation}

\noindent Altogether there are $2N$ points on the grid, with a distance of $h/2$ between them.

\noindent Unlike the standard approach, we note here that the boundary points do not coincide with any grid nodes.

\noindent We consider the following approximation:

\begin{eqnarray}\label{eq2}
\vu_{x} &= & \frac{1}{6h}\left[\left(
                                \begin{array}{ccccccccc}
                                     \ddots &\ddots & \ddots & \ddots &  \ddots &   &   &      \\
                                    & 1  & -8 & {\bf 0} & 8 & -1  &   &       \\
                                    &   &  1 & -8 & {\bf 0} & 8 & -1  &      \\

                                    &   &   & \ddots  & \ddots & \ddots & \ddots  & \ddots  \\
                                \end{array}
                              \right) \right .
       					+  c_1 \left(
                                \begin{array}{ccccccccc}
                                  \ddots &   \ddots & \ddots & \ddots &  \ddots &  \ddots &   &      \\
                                   &  1  & -4 & {\bf 6} & -4 & 1  & 0  &       \\
                                   &  -1  & 4 & -6 & {\bf 4} & -1 & 0 &    \\
                                 &      & \ddots  &  \ddots & \ddots & \ddots & \ddots & \ddots   \\
                                \end{array}
                              \right)     \nonumber   \\
       &&  \hspace{4em}  + c_2 \left . \left(
                                \begin{array}{ccccccccc}
                                  \ddots &   \ddots & \ddots & \ddots &  \ddots &  \ddots &   &      \\
                                   & 0 &  1  & {\bf -4} & 6 & -4  & 1  &       \\
                                   & 0 &  -1  & 4 &  {\bf -6}  & 4 & -1 &    \\
                                 &      & \ddots  &  \ddots & \ddots & \ddots & \ddots & \ddots   \\
                                \end{array}
                              \right)   
                              \right ] \vu  +\vT_e  \;=\; Q \vu  +\vT_e    
\end{eqnarray}

%
%
%
%
%

The truncation error of this scheme is:
\begin{eqnarray} \label{truncation_10}
{ T_e}_{j-\frac{1}{4} } &=&  \dfrac{1}{11520}\left[120(c_{1}+c_{2})h^3u^{(4)}+(-24+60c_{2})h^4u^{(5)}\right]_{x_{j-1/4}}+O(h^{5}) \nonumber \\
\\
{ T_e}_{j+\frac{1}{4} }  &=& \dfrac{1}{11520}\left[-120(c_{1}+c_{2})h^3u^{(4)}+(-24+60c_{1})h^4u^{(5)}\right]_{x_{j+1/4}}+O(h^{5}) \nonumber
\end{eqnarray}
Formally, this is a third-order scheme. Note, however, that the leading term is highly oscillatory.

\subsection{Stability and Accuracy}


\noindent  Unlike the standard finite difference schemes, we need a more refined tool to analyze the BFD methods than estimating the truncation error and applying the Lax-Richtmyer equivalence theorem. Furthermore, since different stencils are used for the points $x_{j-\frac{1}{4} } $ and $x_{j+\frac{1}{4} } $, we cannot use the Fourier analysis directly. Therefore, we use the method proposed in  \cite{ditkowski2015high}, \cite{ditkowski2020error} and \cite{leblanc2020error}.

\bigskip

\noindent The analysis is based on the observation that if  we split the Fourier spectrum into low and high frequencies as follows, let $\omega \in \{-N/2,\ldots,N/2$\}  and:
\begin{equation}\label{} 
\nu \,=\,   \left\{ \begin{array}{lcl}
                       \omega -N, \hspace{1cm}&\omega >0\\
                       \omega+N, \hspace{1cm} &\omega\leq 0
                     \end{array}
\right.
\end{equation}
\noindent Then, for $\omega>0$, we get the following relations:
\begin{equation}\label{S.31} 
 {e}^{i \omega x_{j-1/4}} = i{e}^{i \nu x_{j-1/4}}\;\; {and } \;\; {e}^{i \omega
x_{j+1/4}} = -i{e}^{i \nu x_{j+1/4}} \;.
\end{equation}
\noindent Similarly, for $\omega\leq 0$, we obtain:
\begin{equation}\label{} \nonumber
 {e}^{i \omega x_{j-1/4}} = -i{e}^{i \nu x_{j-1/4}}\;\; { and } \;\; {e}^{i \omega
x_{j+1/4}} = i{e}^{i \nu x_{j+1/4}} \;.
\end{equation}

\noindent Let us assume first that $\omega>0$.  We denote the vectors ${e}^{i \omega \bf{x}}$ and ${e}^{i \nu \bf{x}}$ by:
\begin{equation}\label{} \nonumber
{e}^{i \omega \vx} =  \left(
                          \begin{array}{c}
                            \vdots \\
                            {e}^{i \omega \vx_{j-{1/4}}} \\
                            {e}^{i \omega \vx_{j+1/4}} \\
                            \vdots
                          \end{array}
                        \right) \mbox{, }{e}^{i \nu \vx}= \left(
                          \begin{array}{c}
                            \vdots \\
                            {e}^{i \nu x_{j-{1/4}}} \\
                            {e}^{i \nu x_{j+1/4}} \\
                            \vdots
                          \end{array}
                        \right)
\end{equation}
\noindent We look for eigenvectors in the form of:
\begin{equation}\label{eigenvalue_10} 
\psi_k(\omega) ={\alpha_k}   \frac{   {e}^{i \omega \vx}  }{\sqrt{2 \pi}}    +
{\beta_k}    \frac{  {e}^{i \nu \vx} }{\sqrt{2 \pi}} 
  \end{equation}
\noindent where, for normalization, it is required that
$|\alpha_k|^2+|\beta_k|^2=1 $, $k=1,2$. 
It follows that:
\begin{eqnarray}
Q{e}^{i \omega \vx} ={\rm diag}({\mu}_{1},{\mu}_{2},...,{\mu}_{1},{\mu}_{2}){e}^{i \omega \vx}  \nonumber \\
\\
Q{e}^{i \nu\vx} ={\rm diag}({\sigma}_{1},{\sigma}_{2},...,{\sigma}_{1},{\sigma}_{2}){e}^{i \nu \vx}
 \nonumber 
\end{eqnarray}
%
\noindent where the matrix $Q$ was defined in \eqref{eq2} and:
\begin{eqnarray}
\mu_1  &=& \frac{1}{3h}\left [ -i \left(8 \sin \left(\frac{h  \omega}{2}\right)-\sin (h  \omega) \right)     				-8 \sin ^4\left(\frac{h \omega }{4}\right) 
						\left({c_1}+{c_2} e^{\frac{i h \omega   }{2}}\right)
				 \right ]  \nonumber \\ 
				 \nonumber \\
\mu_2  &=& \frac{1}{3h}\left [  -i \left(8 \sin \left(\frac{h  \omega}{2}\right)-\sin (h  \omega) \right)    
				+8 \sin ^4\left(\frac{h \omega }{4}\right) 
						\left({c_1}e^{-\frac{i h \omega   }{2}}+{c_2} \right)
				 \right ]  \nonumber \\
				 \nonumber \\
\sigma_1  &=& \frac{1}{3h}\left [  i \left(8 \sin \left(\frac{h \omega }{2}\right)+\sin (h \omega )\right) 
					-8 \cos ^4\left(\frac{h \omega }{4}\right) 
							\left({c_1}-{c_2} e^{\frac{i h \omega }{2}}\right)
					 \right ]  \nonumber \\
					 \nonumber \\
\sigma_2  &=& \frac{1}{3h}\left [i \left(8 \sin \left(\frac{h \omega }{2}\right)+\sin (h \omega )\right)
					-8 \cos ^4\left(\frac{h \omega }{4}\right) 
							\left({c_1}e^{-\frac{i h \omega   }{2}}-{c_2} \right)
   					\right ]  \nonumber 
\end{eqnarray}
\noindent We note here that this diagonal matrix is the same under the assumption that $\omega\leq 0$.

\noindent In order to find the coefficients ${\alpha}_{k}$,${\beta}_{k}$ and the  eigenvalues (symbols) $\hat{Q}_k$ for $\omega\neq 0$, we substitute $\psi_k$, defined at \eqref{eigenvalue_10}, into the equation $Q \psi_k = {\hat{Q}}_k  \psi_k $ to get
\begin{equation}
\label{S1.4}
\left. 
\begin{array}{ll}
{\mu}_{1}\frac{\alpha_k}{\sqrt{2 \pi}} {e}^{i \omega {x_{j-\frac{1}{4}}}} +{\sigma}_{1}\frac{\beta_k}{\sqrt{2 \pi}} {e}^{i \nu {x_{j-\frac{1}{4}}}}=
{\hat{Q}}_k\left( \frac{\alpha_k}{\sqrt{2 \pi}} {e}^{i \omega {x_{j-\frac{1}{4}}}} +\frac{\beta_k}{\sqrt{2 \pi}} {e}^{i \nu {x_{j-\frac{1}{4}}}}\right)\\
{\mu}_{2}\frac{\alpha_k}{\sqrt{2 \pi}} {e}^{i \omega {x_{j+\frac{1}{4}}}} +{\sigma}_{2}\frac{\beta_k}{\sqrt{2 \pi}} {e}^{i \nu {x_{j+\frac{1}{4}}}}=
{\hat{Q}}_k\left( \frac{\alpha_k}{\sqrt{2 \pi}} {e}^{i \omega {x_{j+\frac{1}{4}}}} +\frac{\beta_k}{\sqrt{2 \pi}} {e}^{i \nu {x_{j+\frac{1}{4}}}}\right)
\end{array}
  \right.
\end{equation}
\noindent We can simplify \eqref{S1.4}, using the definition $r_{k}=i\dfrac{{\beta}_{k}}{{\alpha}_{k}}$ and obtain the system:
\begin{equation}
\label{S1.5}
\left. 
\begin{array}{ll}
{\mu}_{1}-{\sigma}_{1}r_{k}=
{\hat{Q}}_k\left( 1-r_{k}\right)\\
{\mu}_{2}+{\sigma}_{2}r_{k}=
{\hat{Q}}_k\left( 1+r_{k}\right)\\
\end{array}
  \right.
\end{equation}
\noindent The eigenvalues ${\hat{Q}}_k$, for all $\omega\neq 0$ are:
%
%
\begin{equation}\label{S1.6}
\left. 
\begin{array}{ll}
\hat{Q}_1(\omega) =\dfrac{1}{12h}\bigg(\Omega -\sqrt{2\Delta}\bigg)\\
\hat{Q}_2(\omega) =\dfrac{1}{12h}\bigg(\Omega +\sqrt{2\Delta}\bigg)\\
\end{array}
  \right.
\end{equation}
where
\begin{equation}\label{} \nonumber
\Omega =-4 i (c_1+c_2+1) \sin ( \omega h)+6 (c_2-c_1) \cos (\omega h )+10 (c_2-c_1)\\
\end{equation}
and
\begin{equation}\label{} \nonumber
\left. 
\begin{array}{rcl}
\Delta &=& -80c_1-256+55c_{1}^2+55c_{2}^2-2 (63c_{1}+40)c_{2}+\\
	&& 4 \bigg(15c_{1}^2+c_{1} (16-30 c_{2})+c_{2} (15c_{2}+16)
			+64\bigg) \cos ( \omega h) +\\
	&& \bigg(13c_{1}^2+2c_{1} (8-5c_{2})+c_{2} (13c_{2}+16)\bigg) \cos (2  \omega h)+\\
	&& 8 i (-c_{1}
+c_{2}) \sin (\omega h) \bigg((3c_{1}+3c_{2}+4) \cos (\omega h)+5c_{1}+5c_{2}+28\bigg)
\end{array}
  \right.
\end{equation}
 Using the normalization condition $|\alpha_k|^2+|\beta_k|^2=1 $, $k=1,2$, we choose the coefficients ${\alpha}_{k}$ and ${\beta}_{k}$ to be:
\begin{equation}\label{} 
\alpha_{1}=\dfrac{1}{\sqrt{1+{\vert r_{1}\vert}^2}}\mbox{, } \hspace{2em}
		\beta_{1}=-i\dfrac{r_1}{\sqrt{1+{\vert r_{1}\vert}^2}}\\
\end{equation}
\begin{equation}\label{} 
\alpha_{2}=i \dfrac{\vert r_{2}\vert /r_2}{\sqrt{1+{\vert r_{2}\vert}^2}}\mbox{,  } \hspace{2em}
		\beta_{2}=\dfrac{\vert r_2\vert}{\sqrt{1+{\vert r_{2}\vert}^2}}
\end{equation}

\noindent When $c_1<c_2$, the real part of  $\hat{Q}_1(\omega)$ is strictly positive. Therefore, the scheme is unstable; hence, we will consider the cases where $c_1\geq c_2$.

\subsubsection{Stability proof using eigenvalues and eigenvectors properties}

For any $c_1\geq c_2$, it was numerically shown that the real part of the eigenvalues $\hat{Q}_1(\omega),\hat{Q}_2(\omega)$ is non-positive. Therefore, the scheme is Von-Neumann stable. We next prove that the operator $Q$ has a full set of eigenvectors for any $\omega h$. \\

In the case where $\omega=0$, the corresponding eigenvectors are $\psi_{1}(0)={e}^{i 0\vx},\psi_{2}(0)={e}^{i (N+1)\vx}$ and their corresponding eigenvalues $\hat{Q}_1(0)=0,\hat{Q}_2(0)=(8 (-c1 + c2))/(3 h)$;

Otherwise, by construction, $\psi_{k}(\omega_1)$ and $\psi_{l}(\omega_2)$,  $k,l=1,2$ are orthogonal, for any $\omega_1 \neq \omega_2$.

Unlike the Fourier mode ${e}^{i \omega \vx}$, $\psi_{k}(\omega)$ are not orthogonal, since $\left\langle \psi_{1}(\omega),\psi_{2}(\omega)\right\rangle_h =\sum_{j}\dfrac{h}{2}\overline{\psi_{1}(\omega)_{j} } {\psi_{2}(\omega)_{j}} \neq 0$ where the inner product $\left\langle \cdot,\cdot\right\rangle_h$ is the scalar inner product normalized by $h/2$.

Let us define the cosine of the angle $\theta(\omega)$ between the complex vectors $\psi_{1}(\omega)$ and $\psi_{2}(\omega)$ as follows :

\begin{equation}\label{} \nonumber
\cos(\theta(\omega))=\dfrac{\left |  {\left\langle \psi_{1}(\omega),\psi_{2}(\omega)\right\rangle}_h \right |}{\Vert \psi_{1}(\omega) \Vert \Vert \psi_{2}(\omega)\Vert}
\end{equation}
In fact, we can compute $\cos (\theta(\omega))$ between the vectors using the definition of $r_{1}$ and $r_{2}$, the normalization condition on $\alpha_{k},\beta_{k}$ and \eqref{S.31}.

It yields
\begin{equation}\label{} \nonumber
\cos(\theta(\omega))=
\left | \overline{\alpha_{1}}  {\alpha_{2}}+\overline{\beta_{1}} {\beta_{2}}\   \right | =
\dfrac{1}{\sqrt{1+{\vert r_{1}\vert}^2}}\dfrac{\vert r_{2}\vert}{\sqrt{1+{\vert r_{2}\vert}^2}}\left|  \dfrac{1} {{r_{2}}}+r_1\right|
\end{equation}

It can be shown that $ \left | \cos(\theta(\omega)) \right | <0.4$ for any value $-1\leq c_2\leq c_1\leq 1$ and all $-\pi\leq \omega h \leq \pi$. Therefore, we conclude that the scheme is stable for these values of $c_1$ and $c_2$.  

In the next part of the section, we proceed with the asymptotic analysis of the eigenvalues. For the accuracy analysis, it is sufficient to demonstrate the behaviour of the approximate solution for $\omega h \ll 1$, since, for a solution $u$ which lays in $C^{p}$, its Fourier coefficient $\hat{u}(\omega)$ decays as $const/{\vert \omega \vert ^{p+1}}$. 

\subsection{The $c_1>c_2$ case}

 For $ \omega h \ll 1$, the eigenvalues $\hat{Q}_k(\omega)$, see \eqref{S1.6},  are:
\begin{equation}\label{} \nonumber
\hat{Q}_1(\omega) =-i \omega+\frac{1}{480} i h^4 \omega ^5-\frac{h^5\omega ^6 (c_1+c_2)}{384 (c_1-c_2)} +O(h^6)\\
\end{equation}
and
\begin{equation}\label{}\nonumber
\left. 
\begin{array}{ll}
\hat{Q}_2(\omega) =-\frac{8 (c_1-c_2)}{3 h}+\frac{1}{3} i \omega  (2c_1+2c_2+5)+\frac{1}{2} h \omega ^2 (c_1-c_2)-\frac{1}{9} i h^2 \omega ^3 (c_1+c_2+1)+\\
\\
 \hspace{5em} \frac{1}{24} h^3 \omega ^4 (c_2-c_1)+O(h^4)\\
\end{array}
  \right.
\end{equation}

\noindent The different values for $\alpha_{1},\alpha_{2},\beta_{1}$ and $\beta_{2}$ are

\begin{eqnarray} \nonumber
\alpha_{1}  &=&   1-\frac{(h\omega )^8 (c_1+c_2)^2}{131072 (c_1-c_2)^2}+O\left(h^9\right)\\
 \nonumber \\
\beta_{1}   &=&  -\frac{i(h \omega )^4 (c_{1}+c_{2})}{256 (c_{1}-c_{2})}+\frac{(h \omega )^5 (c_{1}c_{2}+c_{1}+c_{2})}{256 (c_{1}-c_{2})^2}+O(h^6)  \nonumber \\
 \nonumber \\
 \alpha_{2}  &=&  \frac{h\omega  }{4}+
 		\frac{i (h \omega ) ^2}{4c_{1}-4c_{2}} - 
 		\frac{(h \omega )^3 \left(c_{1}^2-2c_{1} (c_{2}-12)+c_{2} (c_{2}+24)+96\right)}{384 (c_{1}-c_{2})^2}+O(h^4)  \nonumber \\
  \nonumber \\
  \beta_{2}  &=&  1 - \frac{1}{32}(h \omega )^2+
  		\frac{(h \omega )^4 \left(-2c_{1}c_{2}+c_{1} (c_{1}+96)+c_{2}^2+96 (c_{2}+2)\right)}{6144 (c_{1}-c_{2})^2}+O(h^6)  \nonumber 
\end{eqnarray}
%
%
%
We note here that it is assumed that the initial data lies only in the low frequencies. \\

\noindent Assuming the initial condition $u(x,0) = {e}^{i \omega x}$, or
\begin{equation}\label{} \nonumber
v_{j-1/4}(0) = {e}^{i \omega x_{j-1/4}}\;, v_{j+1/4}(0) = {
e}^{i \omega x_{j+1/4}} \;\; ; \; \;\;\; h \omega \ll 1
\end{equation}


we obtain, using the expansions for $\hat{Q}_k(\omega),\alpha_{1},\alpha_{2},\beta_{1}$ and $\beta_{2}$,  the expression of the scheme exact solution at $x_{j-1/4}$.
\begin{equation}\label{S1.7}
\left. 
\begin{array}{ll}
{v}_{j-\frac{1}{4}}(t) =
\bigg({e}^{-i \omega t}\bigg[1+\frac{i\omega ^5 t h^4 }{480}-\frac{(\omega h)^5 (c_{1}+c_{2}) (8 t \omega +3 i)}{3072 (c_{1}-c_{2})}\bigg]+O\left(h^6\right)\bigg) {e}^{i \omega x_{j-\frac{1}{4}}}+\\
\\
\hspace{5em}
{e}^{-i  \omega t}\left(\frac{(\omega h)^5 (c_{1}c_{2}+c_{1}+c_{2}) }{256 (c_{1}-c_{2})^2}-\frac{i(\omega h)^4 (c_{1}+c_{2}) }{256 (c_{1}-c_{2})}+O\left(h^6\right)\right) {e}^{i \nu x_{j-\frac{1}{4}}}
\end{array}
  \right.
\end{equation}
A similar expression holds for $x_{j+1/4}$. Therefore, the scheme has a fourth-order convergence rate, for all $c_1>c_2$.

We note here that the order of magnitude of $\alpha_1$ and $\beta_2$ is 1, whereas $\beta_{1}$ is of order $h^4$ and $\alpha_{2}$ is of order $h$. Hence, since the exact solution is ${e}^{-i \omega t}$, the exact solution of the equation is well approximated by the low frequencies represented by $\omega$, while the truncation error lies with the high frequencies $\nu$.

\subsection{The  $c_1=c_2$ case}

We now explore the case $c_1=c_2=c$. In this case, the $\hat{Q}_k(\omega)$ become:
\begin{equation}\label{} \nonumber
\left. 
\begin{array}{ll}
\hat{Q}_1(\omega) =\dfrac{1}{3h}\bigg(-2 \sqrt{2} \sin \left(\frac{h \omega }{2}\right) \sqrt{-(c+2) (c \cos (h \omega )+c+4)}+i (2 c+1) \sin (h \omega )\bigg)\\
\\
\hat{Q}_2(\omega) =\dfrac{1}{3h}\bigg(2 \sqrt{2} \sin \left(\frac{h \omega }{2}\right) \sqrt{-(c+2) (c \cos (h \omega )+c+4)}+i (2 c+1) \sin (h \omega )\bigg)\\
\end{array}
  \right.
\end{equation}
\noindent Following the same procedure as in the previous section, we conclude that for $\omega h \ll 1$, the eigenvalues are:
\begin{equation}\label{} \nonumber
\hat{Q}_1(\omega) =-i \omega+\frac{i (1-2 c) h^4 \omega ^5}{240 (c+2)}+O(h^6)\\
\end{equation}
\noindent and
\begin{equation}\label{} \nonumber
\left. 
\begin{array}{ll}
\hat{Q}_2(\omega) =\frac{1}{3} i (4 c+5) \omega-\frac{1}{9} i (2 c+1) h^2 \omega ^3+\frac{i \left(8 c^2+26 c+5\right) h^4 \omega ^5}{720 (c+2)}+O(h^6)\\
\end{array}
  \right.
\end{equation}
$\alpha_{1},\alpha_{2},\beta_{1}$ and $\beta_{2}$ are:
\begin{equation} \nonumber
{\alpha}_{1} =  1-\frac{ c^2 \omega ^6 h^6  }{8192 (c+2)^2} + O\left(h^{10}\right),  \quad
{\beta}_{1} = \frac{c \omega ^3 h^3 }{64 (c+2)} + 
		\frac{c (c-2)  \omega ^5  h^5 }{1024 (c+2)^2}+O\left(h^7\right)
\end{equation}
and
\begin{equation} \nonumber
{\alpha}_{2} =  \frac{c  \omega h }{4 c+8} - 
		\frac{c \left(c^2+16 c+40\right) \omega^3 h^3 }{384 (c+2)^3} + O\left(h^5\right) ,  \;
{\beta}_{2 } =  1 - \frac{  c^2 \omega ^2   h^2 }{32 (c+2)^2} + O\left(h^4\right)
\end{equation}


%
%
%
%
\noindent  Assuming the initial condition $u(x,0) = {e}^{i \omega x}$, or
\begin{equation}\label{} \nonumber
v_{j-1/4}(0) = {e}^{i \omega x_{j-1/4}}\;, v_{j+1/4}(0) = {
e}^{i \omega x_{j+1/4}} \;\; ; \; \;\;\; \omega h\ll 1
\end{equation}
%
%
%
%
\noindent and we obtain,  for $\omega>0$,  the expression of the scheme exact solution at $x_{j-1/4}$.
\begin{equation}\label{S1.8} 
\left. 
\begin{array}{ll}
{v}_{j-\frac{1}{4}}(t) =
\left [  {e}^{-i\omega t}  
					{  {e}^ {i \left (  \frac{ (1-2 c)  \omega^5 h^4 }{240 (c+2)} - 
 			\frac{  \left(2 c^2-6 c+1\right)  \omega^7 h^6 }{4032 (c+2)^2}  + 
 					O\left( \omega^9 h^{8}\right)   \right )t }  } \right ] {e}^{i \omega x_{j-\frac{1}{4}}} +
 				 \\
 \\ 				 
 				{   \qquad
		\Bigg[   \Bigg (   \frac{c^2  w^4  h^4 }{256 (c+2)^2}+O\left(\omega^6 h^6\right)   \Bigg ) \times    } 
		  		{ 
			\Bigg (  {e}^{-i \left ( \omega +O\left(\omega^5 h^4\right)   \right )t}    - 
			  {e}^{i \left (  \frac{1}{3} (4 c+5) \omega   + O\left(\omega^3 h^2\right)   \right )t} 
 			\Bigg)    \Bigg]  {e}^{i \omega x_{j-\frac{1}{4}}}  } +
 				 \\
 \\ 				 
		{   \qquad
		\Bigg[   \Bigg (  \frac{c  \omega ^3 h^3}{64 (c+2)} + 
				O\left(\omega^5 h^5\right)  \Bigg  ) 
				 \times  } 
		  		{  
			\Bigg (  {e}^{-i \left ( \omega +O\left(\omega^5 h^4\right)   \right )t}    - 
			  {e}^{i \left (  \frac{1}{3} (4 c+5) \omega   + O\left(\omega^3 h^2\right)   \right )t} 
 			\Bigg)    \Bigg]{ e }^{i \nu x_{j-\frac{1}{4}}}  }  \\
\end{array}
  \right.
\end{equation}
%
Similar expressions holds for $x_{j+1/4}$ and $\omega\le 0$.
There are three sources for the error. The third one is a third-order, bounded in time high-frequency term. This term, however, can be eliminated in a post-processing stage. The second one is a fourth-order, low-frequency one. This term is also bound in time. The first term is a phase error. It grows linearly in time, at least until the error is approximately in order of 1. The phase error is fourth-order; however, it can be made six-order by taking $c=1/2$.

\bigskip

Therefore, for the $c=1/2$ case, what we expect to see for short to moderate times is a third-order error, bounded in time, which can be filtered using a post-processing filter to get a bounded-in-time fourth-order error. When the six-order linearly growing phase error becomes dominant, this scheme behaves like a standard six-order finite difference scheme. It should be noted that the error in the standard finite difference scheme is a linearly growing phase error. 

We demonstrate this behavior in the next section.

\subsection{Numerical examples}

\begin{figure}[htp]
\centering
\begin{center}
\begin{tabular}{lll}
a&b\\
\hspace{-2em}
\includegraphics[clip,scale=0.13]{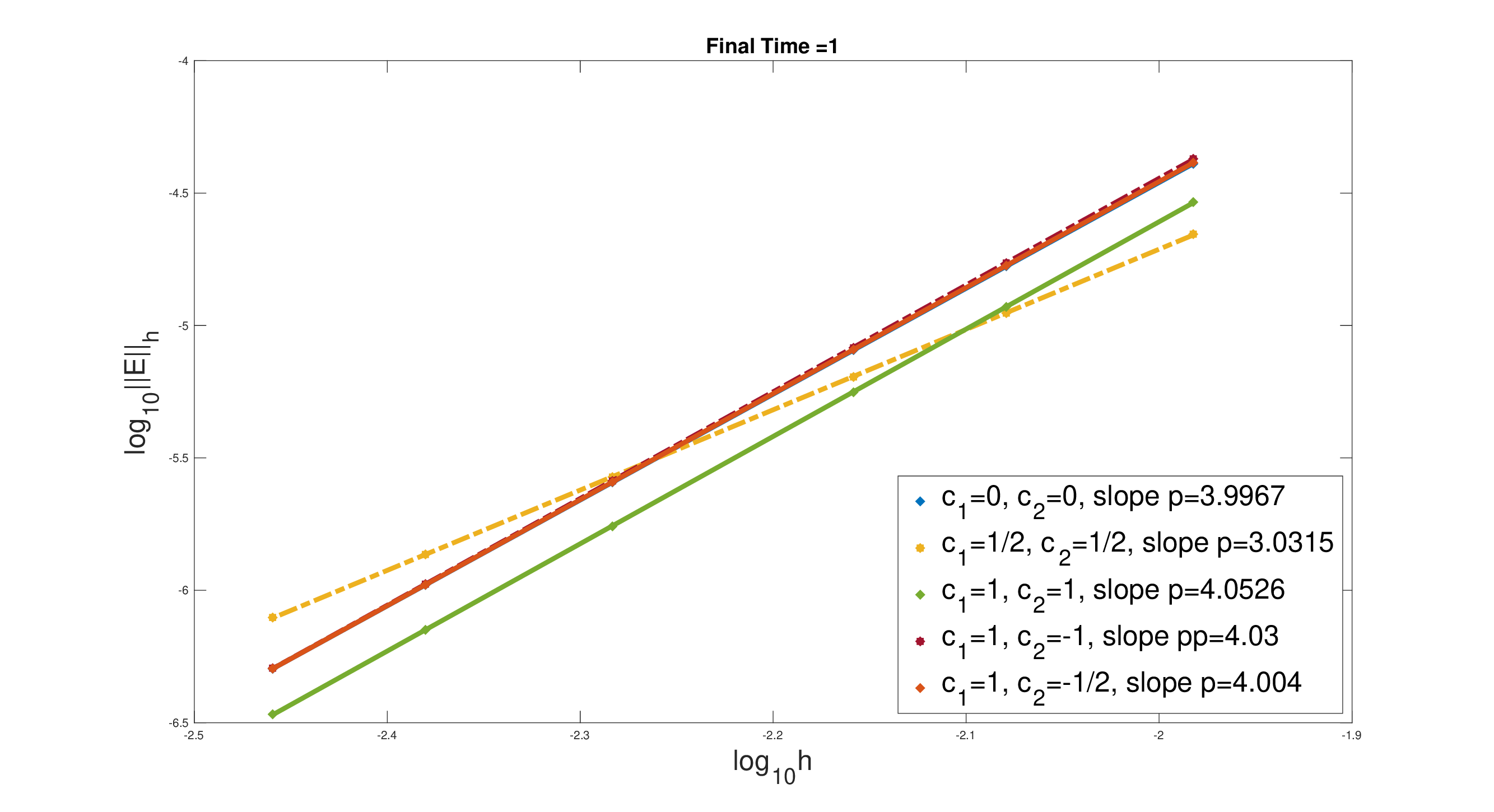}&
\hspace{-2em}
\includegraphics[clip,scale=0.13]{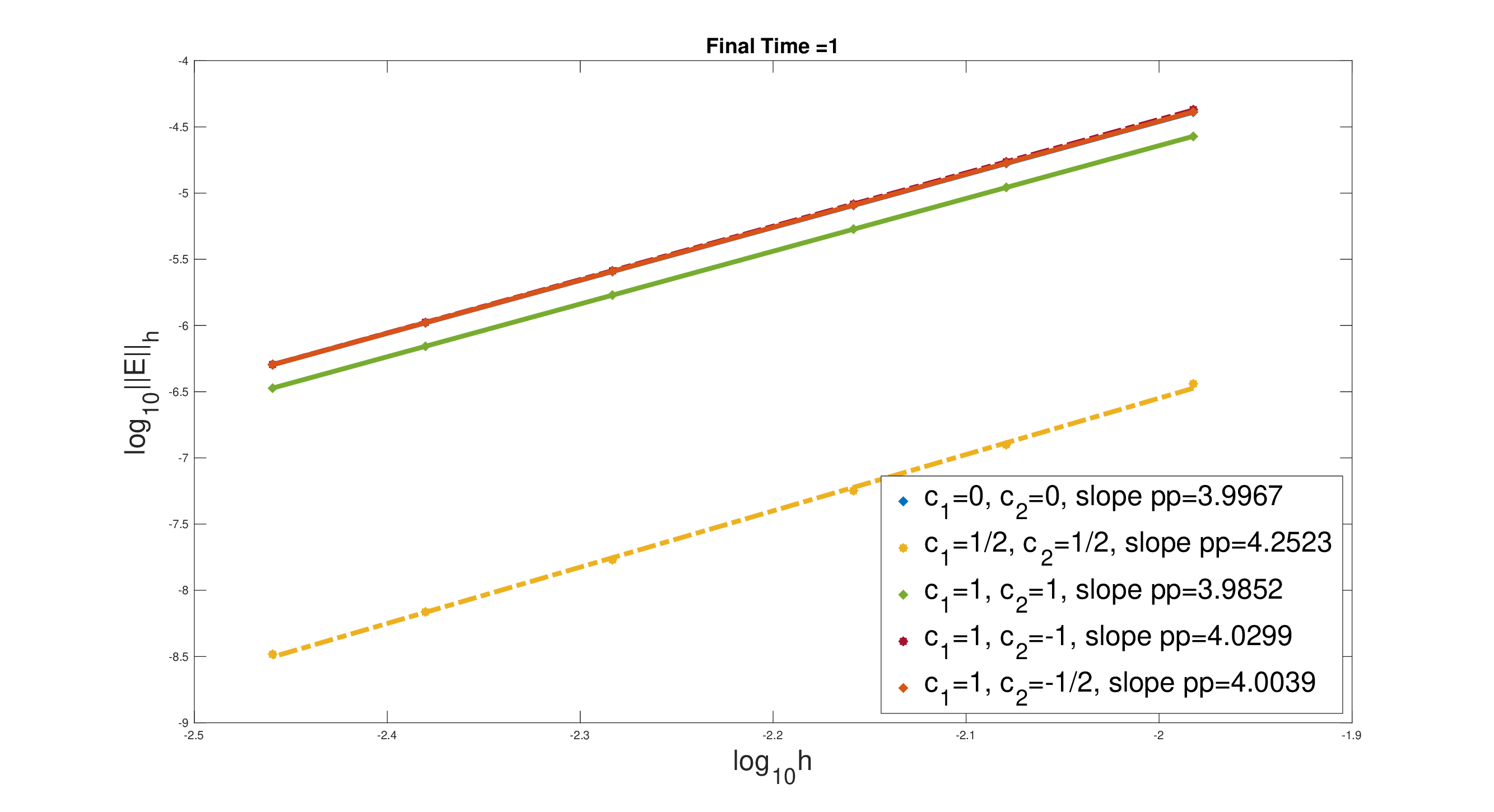}
\end{tabular}
\end{center}
\caption{Scheme \eqref{eq2} Convergence plots, $\log_{10} \|\vE\|$ vs. $\log_{10} h$ ,for  different values of $c_1,c_2$. Final time, $T=1$ -a: no post-processing; b: spectral post-processing}\label{fig:Numerical_examples_10}
\end{figure}

In this section, we illustrate the properties of the BFS method presented in the previous section using a few numerical examples.

\begin{equation}
\label{Numerical_examples_10}
\left\lbrace 
\begin{array}{ll}

\dfrac{\partial u}{\partial t}+ \dfrac{\partial u}{\partial x } =0\mbox{ , }x\in \left( 0,2\pi \right) \mbox{ , }t\geq 0\\
u(x,0)=\exp\cos(2\pi(x))\\

\end{array}
  \right.
\end{equation}
The exact solution to this problem is $u(x,t)=\exp\cos(2\pi(x-t))$. 
Scheme \eqref{eq2} was run, for  different values of $c_1,c_2$, on the interval $[0,1]$ with $N=48,60,72,96,120,144$, and he final time was $T=1$. A six-order explicit Runge-Kutta method was used for time propagation.  The convergence plots are presented in Fig. \ref{fig:Numerical_examples_10}.

As shown in Fig. \ref{fig:Numerical_examples_10}a, where $c_1=-c_2$, the truncation error, \eqref{truncation_10}, and the scheme are fourth-order. In the case $c_1=1$ and $c_2=-1/2$, where $c_1>c_2$, it is a fourth-order scheme, although the truncation error is only third-order. When $c_1=c_2=1/2$, the scheme is indeed third-order; however, when $c_1=c_2=1$, the convergence rate is four due to the cancellation of the leading terms of the bounded parts of the error at integer times, see equation \eqref{S1.8}. After post-processing, see Fig. \ref{fig:Numerical_examples_10}b, the scheme has a fourth-order convergence rate. The error for the $c_1=c_2=1/2$ case is two orders of magnitude smaller since it is caused only by the bounded part of the error, and the linearly growing, six-order phase error is neglectable.

\begin{figure}[htp]
\centering
\begin{center}
\begin{tabular}{lll}
a&b\\
\hspace{-2em}
\includegraphics[clip,scale=0.13]{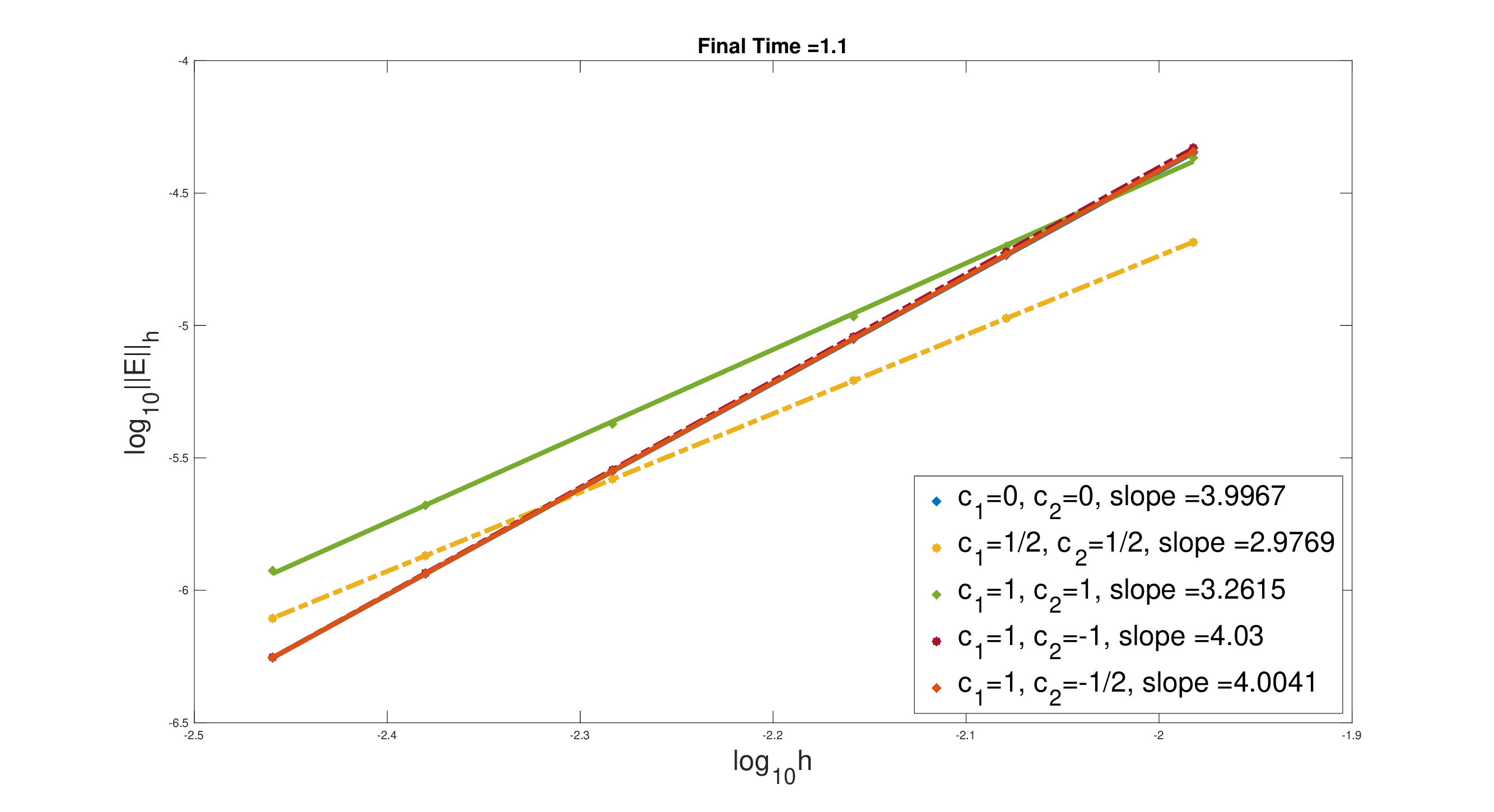}&
\hspace{-2em}
\includegraphics[clip,scale=0.13]{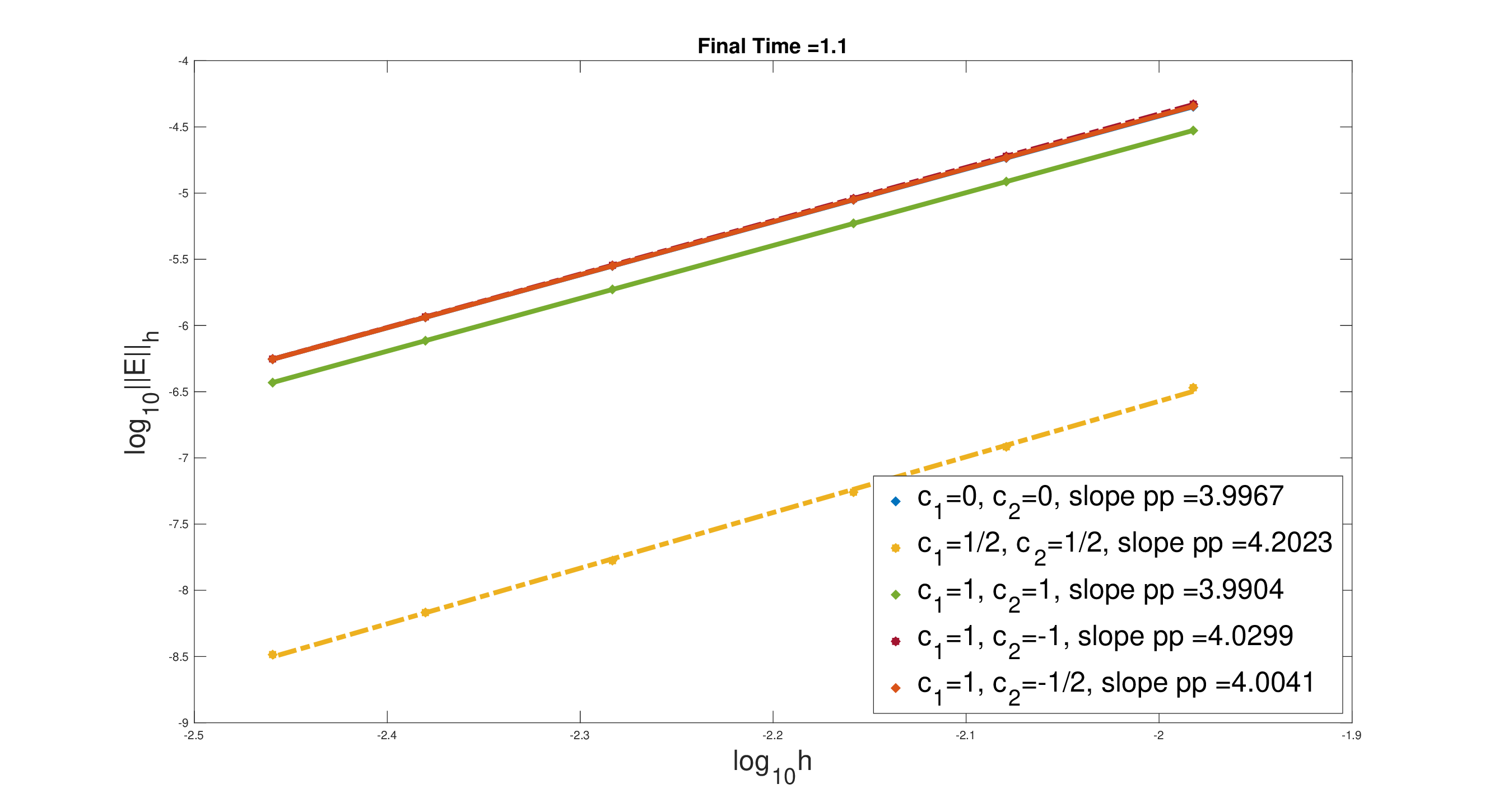}
\end{tabular}
\end{center}
\caption{Scheme \eqref{eq2} Convergence plots, $\log_{10} \|\vE\|$ vs. $\log_{10} h$ ,for  different values of $c_1,c_2$. Final time, $T=1.1$ -a: no post-processing; b: spectral post-processing}\label{fig:Numerical_examples_20}
\end{figure}
When the final time is $T=1.1$, see Fig. \ref{fig:Numerical_examples_20}, the scheme has a third-order convergence rate for $c_1=c_2=1$, which can be filtered to be fourth-order after the post-processing. Since the phase error is of fourth order, the error is comparable to the other $c_1>c_2$ cases.

\begin{figure}[htp]
\centering
\begin{center}
\begin{tabular}{lll}
a&b\\
\hspace{-2em}
\includegraphics[clip,scale=0.14]{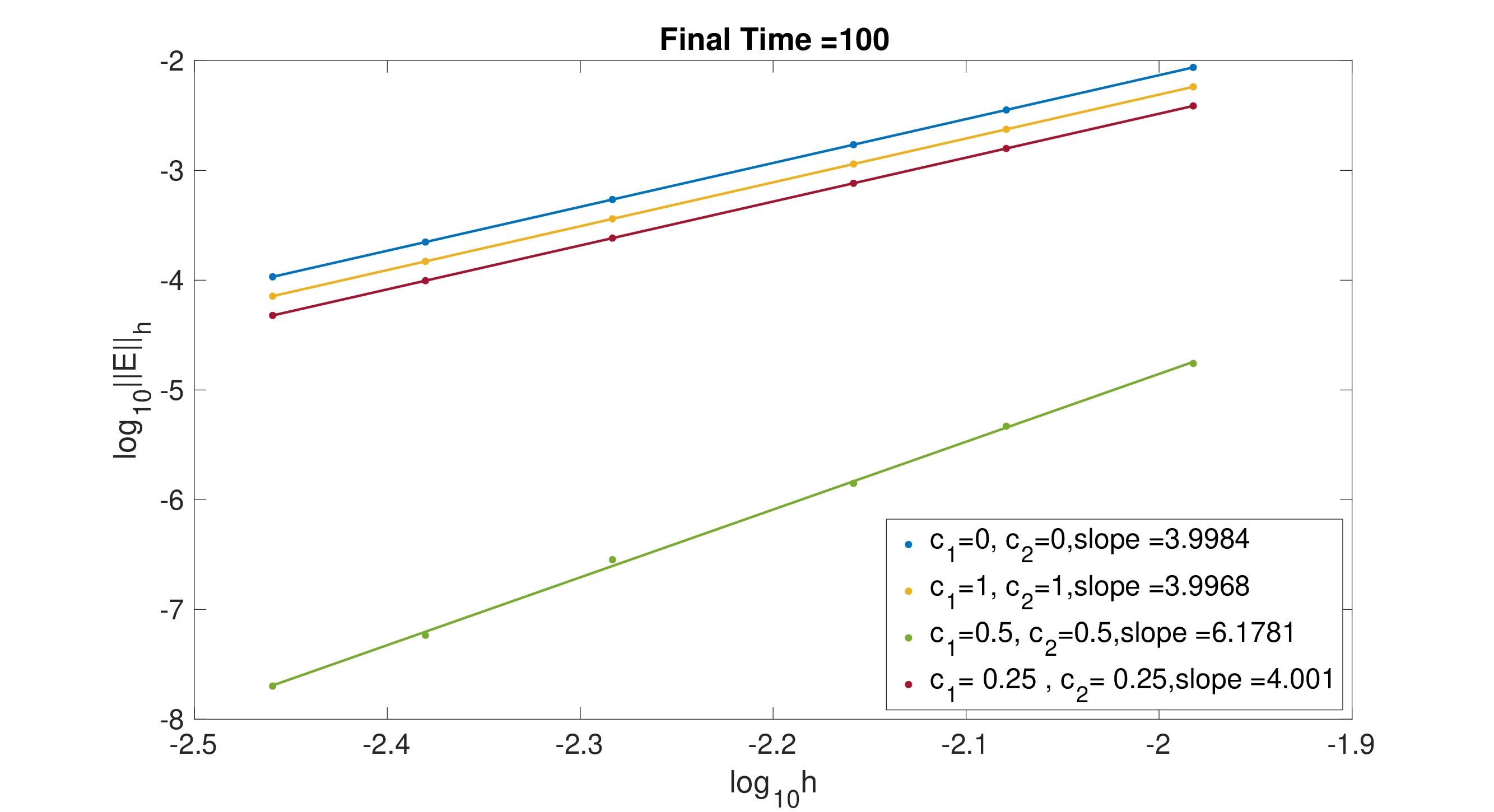}&
\hspace{-2em}
\includegraphics[clip,scale=0.14]{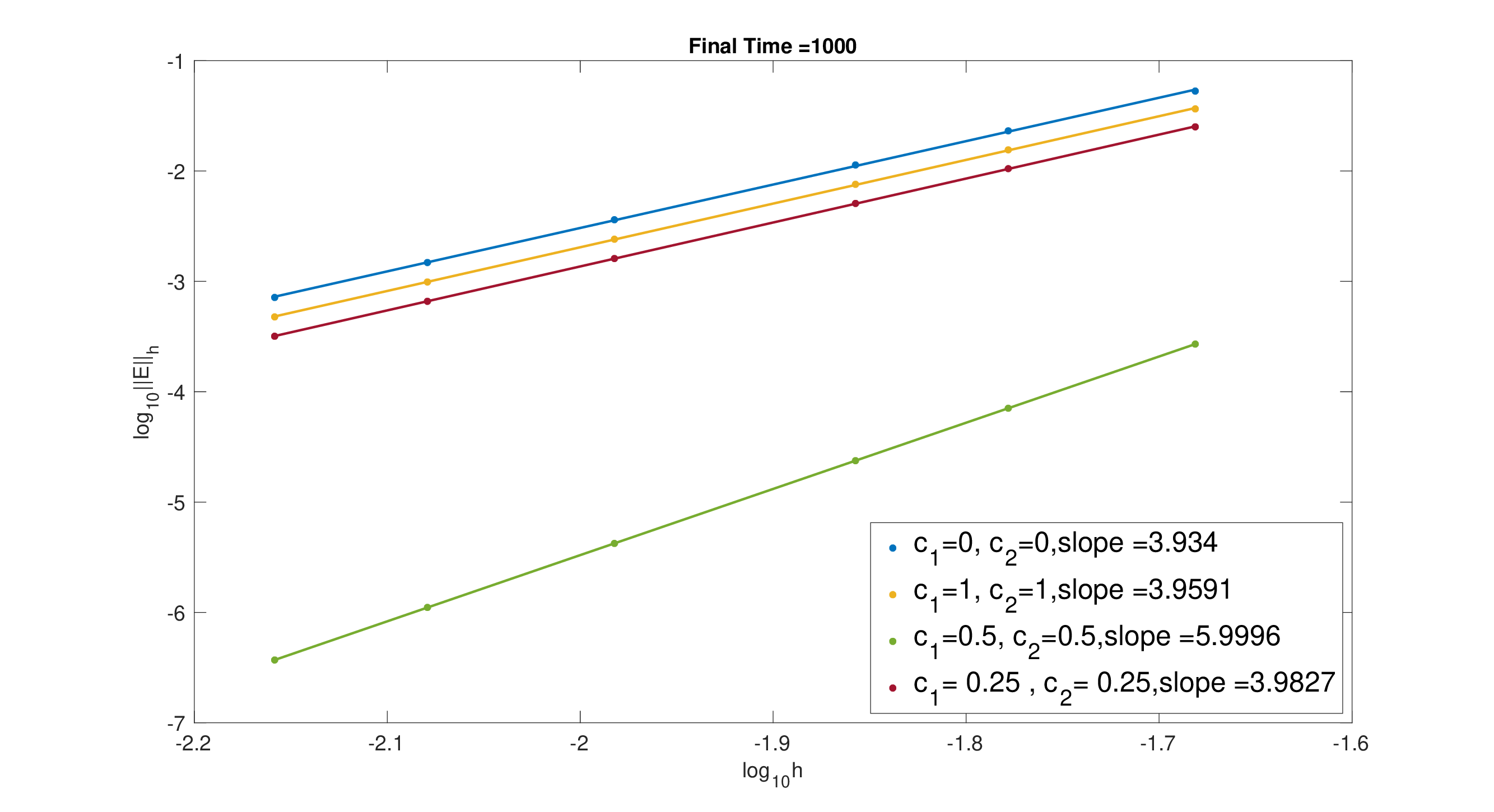}
\end{tabular}
\end{center}
\caption{Scheme \eqref{eq2}, Convergence plots, $\log_{10} \|\vE\|$ vs. $\log_{10} h$ - Periodic BC for $c_{1}=c_{2}$, with post processing, for different final times. a:  $T=100$, b: $T=1000$.}\label{fig:Numerical_examples_30}
\end{figure}

To examine the long-time behavior of the scheme for $c_1=c_2$, we ran the example for $T=100$ and $T=1000$: the analysis, eq. \eqref{S1.8} shows that for $c_1=c_2=1/2$, the scheme convergence rate becomes six-order when the phase error becomes dominant. The convergence rates presented in Fig. \ref{fig:Numerical_examples_30} demonstrate this phenomenon. 

\begin{figure}[htp]
\centering
\begin{center}
\begin{tabular}{lll}
a&b\\
\hspace{-2em}
\includegraphics[clip,scale=0.14]{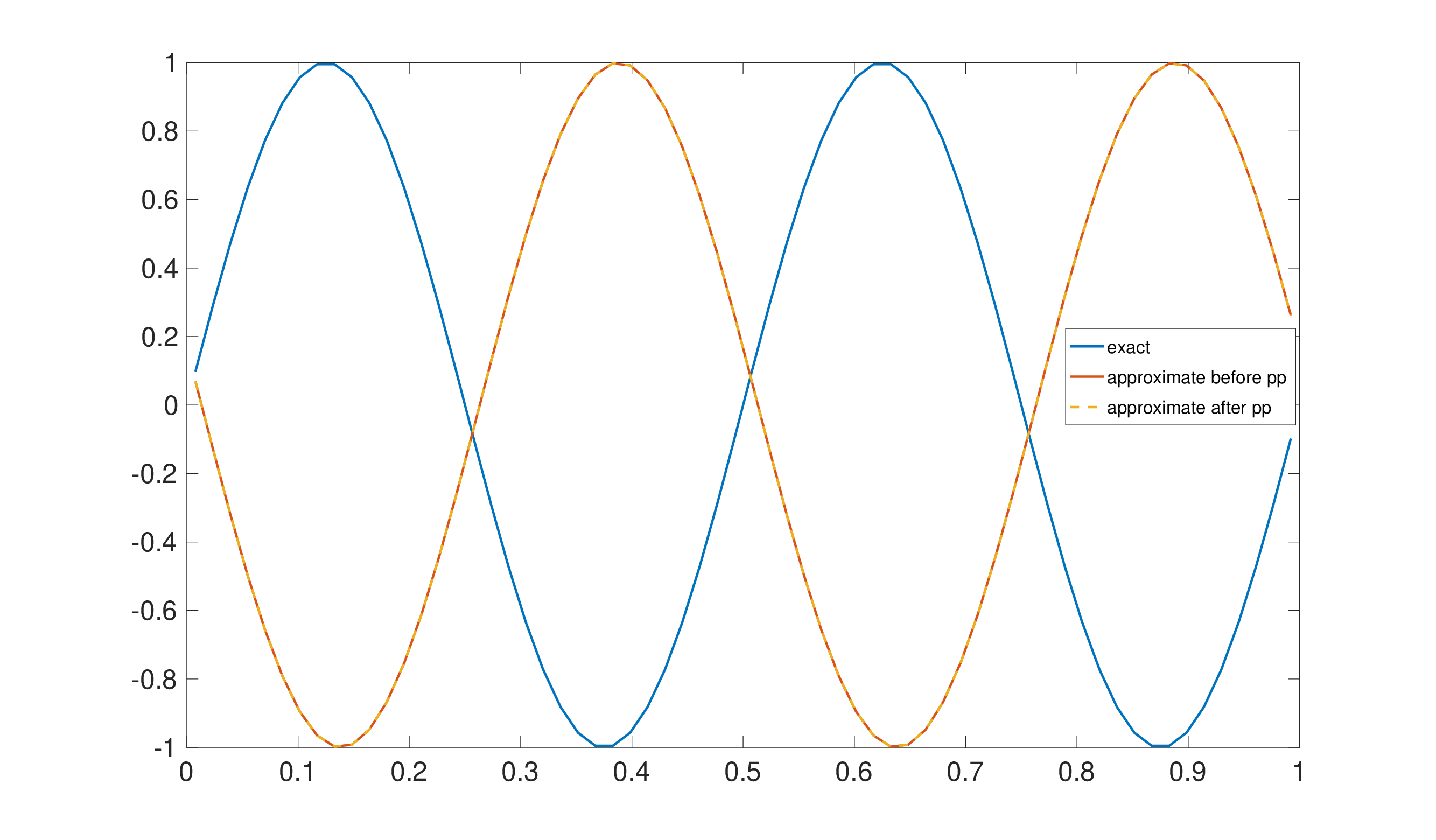}&
\hspace{-2em}
\includegraphics[clip,scale=0.14]{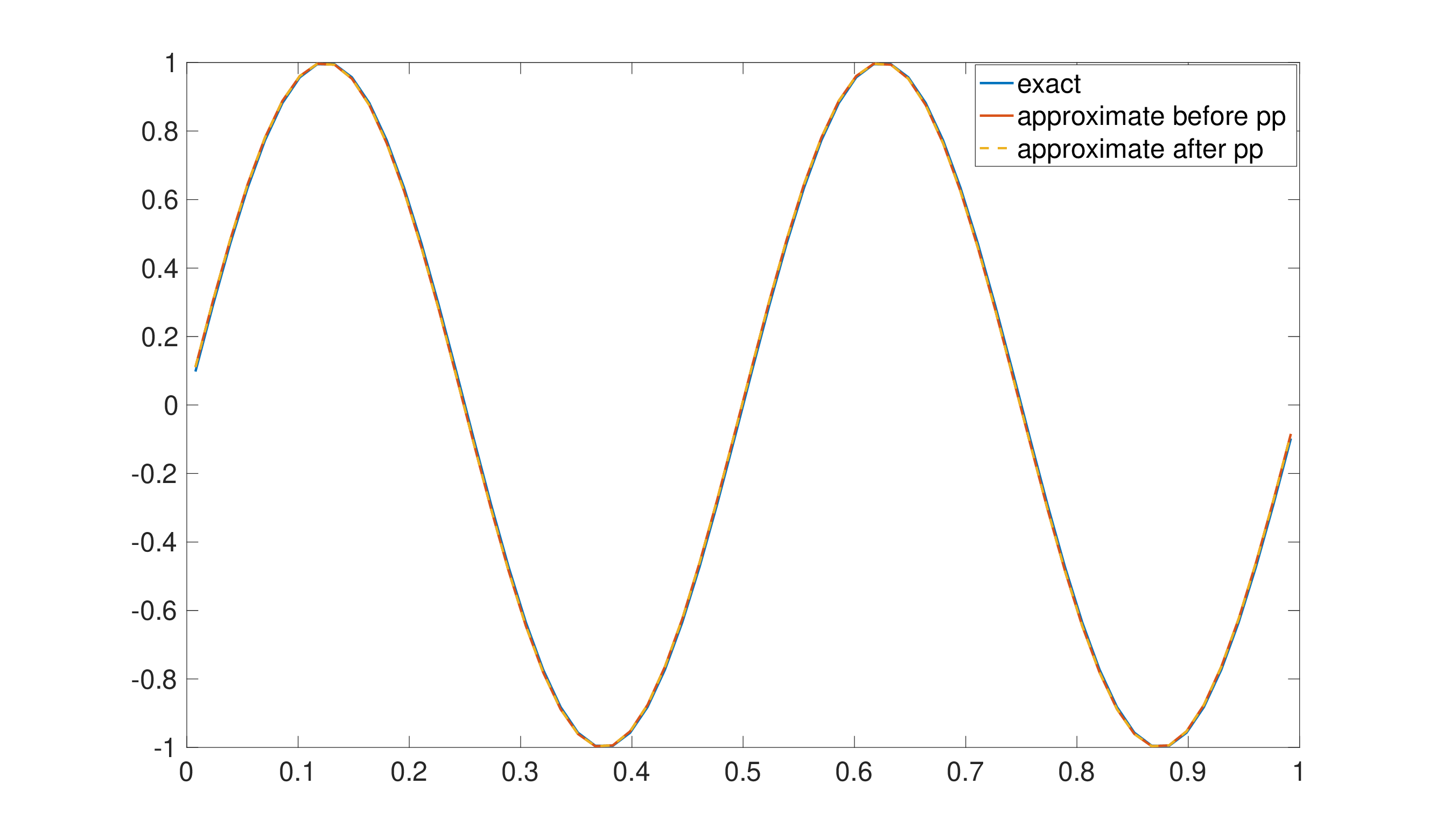}
\end{tabular}
\end{center}
\caption{The exact numerical solutions for $u(x,t)=\sin(4\pi(x-t))$, final time  $T=4800$  and $N=32$ with post processing. a: standard fourth-order scheme, $c_{1}=c_{2}=0$,  b: $c_{1}=c_{2}=1/2$.}\label{fig:Numerical_examples_40}
\end{figure}

To illustrate the effect of the phase error for long-time propagations, we ran the scheme with the initial condition $f(x)=\sin(4\pi(x))$ with $N=32$ (64 degrees of freedom) to  $T=4800$. We compare the cases  $c_1=c_2=0$, the standard fourth-order FD scheme,  and $c_1=c_2=1/2$, with and without post-processing. The results are shown in Fig.  \ref{fig:Numerical_examples_40}. While the standard scheme has a $180^{\circ}$ phase error, the $c_1=c_2=1/2$, solution is undistinguished from the exact solution, at least in the eyeball norm. This observation is valid for both the un-post-processed and the post-processed cases. 

\begin{figure}[htp]
\centering
\begin{center}
\begin{tabular}{lll}
a&b\\
\hspace{-2em}
\includegraphics[clip,scale=0.36]{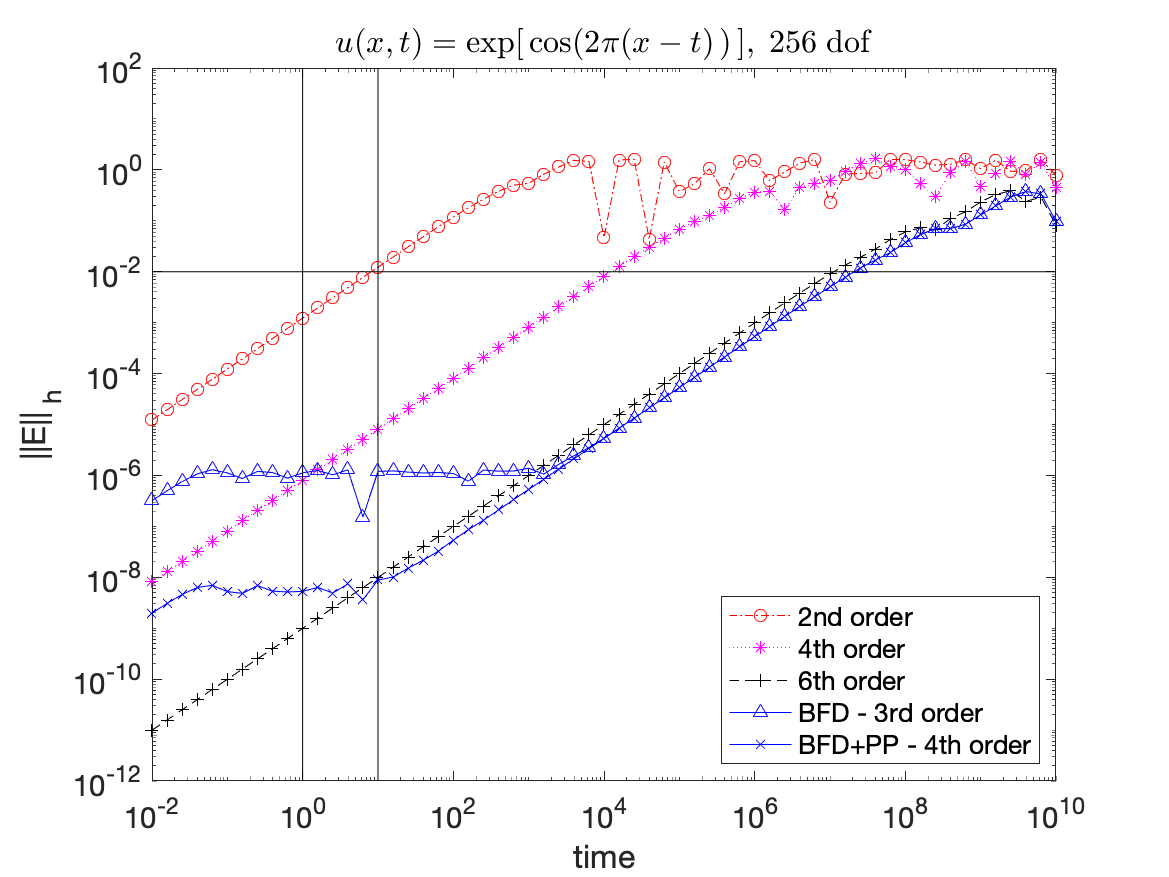}&
\includegraphics[clip,scale=0.36]{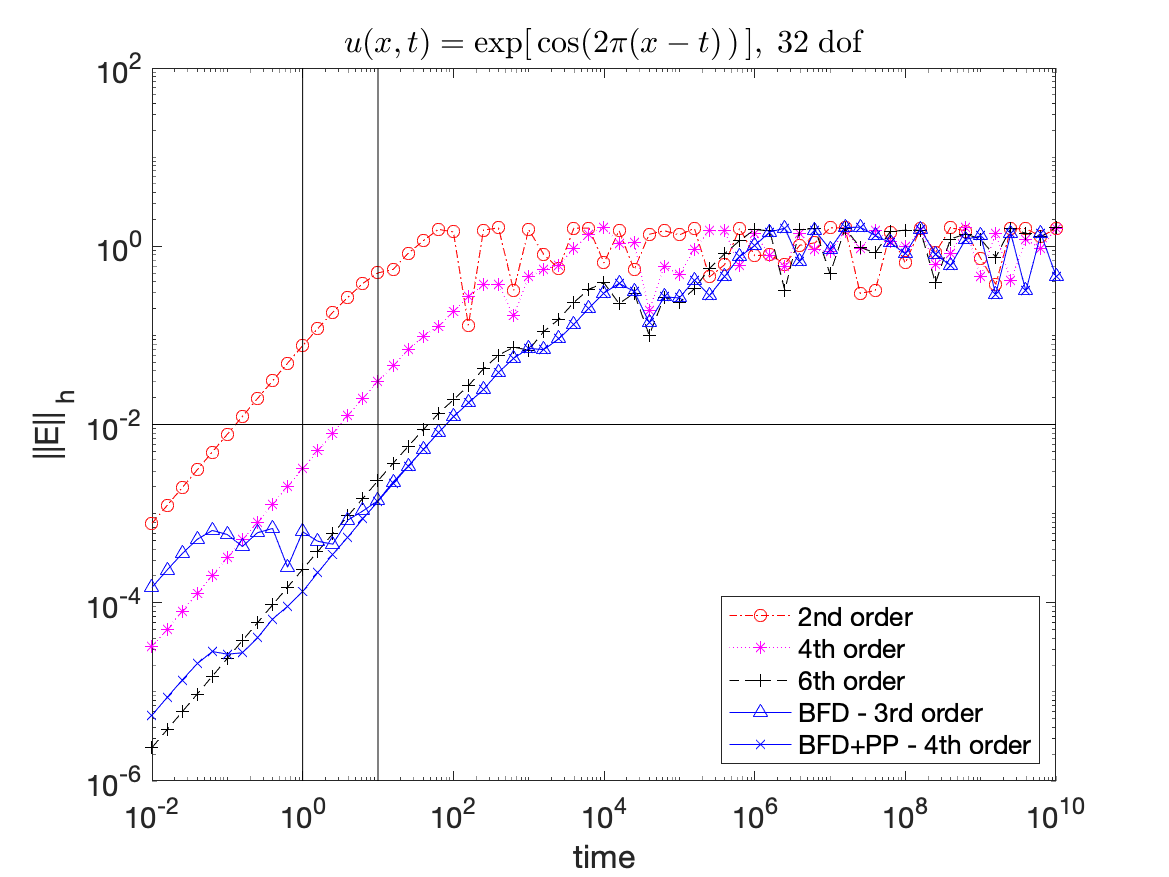}
\end{tabular}
\end{center}
\caption{The evolution of the error in time for the standard second, fourth, and six-order finite difference schemes and the BFD scheme,  \eqref{eq2}
, with  $c_1=c_2=1/2$ with and without post-processing. a: for $N=128$ (256 degrees of freedom), b: for $N=16$ (32 degrees of freedom).}\label{fig:Numerical_examples_50}
\end{figure}
To examine the way the errors evolved in time, we solved the problem \eqref{Numerical_examples_10} to $T=10^{10}$ using the standard second, fourth, and six-order finite difference schemes and the BFD scheme,  \eqref{eq2}, with  $c_1=c_2=1/2$ with and without post-processing. We computed the time propagation by calculating the matrix $\exp(Qt)$. For the standard schemes, the errors grow linearly in time, as they are phase errors. For the  $c_1=c_2=1/2$ case, the third-order error and fourth-order after post-processing are bounded in time. When the six-order linearly growing phase error becomes dominant, this scheme behaves like a standard six-order finite difference scheme. The magnitudes of the errors and the rate at which the phase errors grow in time depend on $N$ (or $h$). Thus, for high resolution, $N=128$ (256 degrees of freedom), the phase error becomes dominant at $t \approx 10$ for the post-processed solution and at $t \approx 1000$  for the unprocessed one. Thus, for low resolution, $N=16$ (32 degrees of freedom), the phase error becomes dominant at $t \approx 1$ for the unprocessed solution and almost immediately for the processed one, as illustrated in Fig.  \ref{fig:Numerical_examples_50}.

\section{Equivalence to DG Method} \label{sec:DG}



This section shows that Scheme \eqref{eq2} is also a nodal-based, $p=1$, Discontinuous Galerkin (DG) scheme.  This DG scheme has a fourth-order (after post-processing) order of convergence for short times and six-order for long times, whereas the standard DG scheme with $p=1$ is typically a second-order method.

\subsection{Standard DG scheme}  \label{Sec:Standard DG}

Discontinuous Galerkin (DG) methods are a particular class of finite element methods using discontinuous basis functions. We start by deriving a standard DG method as illustrated in \cite{zhang2003analysis}  and \cite{cockburnDG}.

Consider the transport equation \eqref{eq1}, with periodic boundary conditions.
\begin{equation}
\label{w_1}
\left\lbrace 
\begin{array}{ll}

\dfrac{\partial u}{\partial t}+ \dfrac{\partial u}{\partial x } =0\mbox{ , }x\in \left( 0,2\pi \right) \mbox{ , }t\geq 0\\
u(x,0)=f(x)\\

\end{array}
  \right.
\end{equation}
We assume the following mesh to cover the computational domain $[0,2\pi]$, consisting of cells
\begin{equation} \nonumber
I_{j}=\left[ x_{j-{1/2}},x_{j+{1/2}}\right] \;,\;\;\;\;\;j=1,...,N
\end{equation}
 where
 \begin{equation} \nonumber
0=x_{1/2}<x_{3/2}<...<x_{N+1/2}=2\pi
\end{equation}
 The center of each cell is located at $x_{j}=\frac{1}{2}\left( x_{j-{1/2}}+x_{j+{1/2}}\right) $ and the size of each cell is $\Delta x_{j}=x_{j+{1/2}}-x_{j-{1/2}}$. We consider a uniform mesh,  hence $h=\Delta x=\frac{2\pi}{N}$. \\

The discontinuous Galerkin scheme is defined as follows:
Find $u,v\in V_{\Delta x}$ (where $V_{\Delta x}=\lbrace v:$ $v$ is a polynomial of degree at most $p$ for $x\in I_{j},j=1,...,N\rbrace$) such that:\\
\begin{equation} \label{w_5}
\left.
\begin{array}{ll}
\int_{I_j}u_{t}v dx-\int_{I_j}uv_{x}dx +
		\hat{u}_{j+\frac{1}{2}} v^{-}_{j+1/2} -
		\hat{u}_{j-\frac{1}{2}}  v^{+}_{j-1/2}=0\\
  \end{array}
  \right.
\end{equation}
for all $v\in V_{\Delta x}$. In our case, the information flows from left to right; we chose the upwind flux:
\begin{equation}
\label{w_4}
\left.
\begin{array}{ll}
\hat{u}_{j+\frac{1}{2}}=u^{-}_{j+\frac{1}{2}} \\
\hat{u}_{j-\frac{1}{2}}=u^{-}_{j-\frac{1}{2}} \\
  \end{array}
  \right.
\end{equation}
 Where a linear element basis, $p=1$,  is used on an equidistant grid, $x_{j-1/4}$ and $x_{j+1/4}$, $u(x)$ and $v(x)$ have the form
\begin{eqnarray} \label{nodal_u_v_def}
u(x) &=& u_{j-1/4} \varphi_{j-1/4}  + u_{j+1/4} \varphi_{j+1/4} \nonumber \\
v(x) &=& v_{j-1/4} \varphi_{j-1/4}  + v_{j+1/4} \varphi_{j+1/4}
\end{eqnarray}
where the Lagrange interpolating polynomials $\varphi_{j \pm 1/4}$ are
\begin{equation}
\label{66}
\left.
\begin{array}{ll}

   \varphi_{j-1/4}= -\frac{2}{\Delta x}(x-x_{j+\frac{1}{4}})\\
    \varphi_{j+1/4}=\frac{2}{\Delta x}(x-x_{j-\frac{1}{4}})\\

  \end{array}
  \right.
\end{equation}
Collecting the coefficients of $v_{j-1/4}$ and $v_{j+1/4}$ yields the equation for the time derivative of $u_{j-1/4}$ and $u_{j+1/4}$
\begin{equation}
\label{w_6}
\left.
\begin{array}{ll}

   \begin{bmatrix}
           u_{j-1/4} \\
           u_{j+1/4} \\
          
         \end{bmatrix}_{t} &=\left( A\begin{bmatrix}
           u_{j-5/4} \\
           u_{j-3/4} \\
          
         \end{bmatrix} +B\begin{bmatrix}
           u_{j-1/4} \\
           u_{j+1/4} \\
          
         \end{bmatrix}+C\begin{bmatrix}
           u_{j+3/4} \\
           u_{j+5/4} \\
          
         \end{bmatrix}
 \right)
  \end{array}
  \right.
\end{equation}
where 
\begin{equation}   \label{eq5}
 A = \frac{1}{4 h}\begin{bmatrix}
          -5 & 15\\
          1 & -3 \\
         \end{bmatrix} \;,\;\;\;\;\; 
B = \frac{1}{4 h}\begin{bmatrix}
           -7 & -3 \\
          11 & -9\\
        \end{bmatrix}   \;,\;\;\;\;\;  
C = \frac{1}{4 h}\begin{bmatrix}
          0 & 0\\
          0 & 0 \\
           \end{bmatrix}   
\end{equation} 
%
We note here that the standard DG scheme has a second order accuracy ($k+1$ order, where $k=1$ in our example).

\subsection{Proof of equivalence}

We can write Scheme \eqref{eq2} in the same form as  \eqref{w_6}, by assembling the matrices $A,B$ and $C$ are defined as follows:
\begin{eqnarray} \label{DG_matrices}
 A &=& \frac{1}{6h}\begin{bmatrix}
           -1-c_{1}&8+4c_{1}-c_{2} \\
           c_{1}&-1-4c_{1}+c_{2} \\
         \end{bmatrix}  \;,\;\;\;\;
B \;=\; \frac{1}{6h}\begin{bmatrix}
           -6c_{1}+4c_{2}&-8+4c_{1}-6c_{2} \\
           8+6c_{1}-4c_{2}&-4c_{1}+6c_{2}\\
        \end{bmatrix}      \nonumber    \\
&& \hspace{7em} C\;=\;  \frac{1}{6h}\begin{bmatrix}
           1-c_{1}+4c_{2}&-c_{2} \\
           -8+c_{1}-4c_{2}&1+c_{2} \\
           \end{bmatrix}        
\end{eqnarray}
\noindent Hence, our goal is to find the corresponding weak formulation of the problem numerical fluxes and possibly other penalty terms, such that the BFD scheme can be viewed as a form of DG scheme.

\noindent As done in Section \ref{Sec:Standard DG}, we choose a linear element basis \eqref{nodal_u_v_def} for the test and trial functions. 

\noindent After replacing boundary terms with fluxes and the test function by its values inside the cell, the scheme becomes:
\begin{equation}
\label{eq6}
\left.
\begin{array}{ll}
\int_{I_j}u_{t}v dx-\int_{I_j}uv_{x}dx+\hat{u}_{j+1/2}v^{-}_{j+1/2}-\hat{u}_{j-1/2}v^{+}_{j-1/2}=0
  \end{array}
  \right.
\end{equation}
\noindent By using the upwind flux and adding in Eq.\eqref{eq6} and adding all the possible penalty terms with general coefficients, we obtain the following scheme:
\begin{equation}
\label{eq7}
\left.
\begin{array}{ll}
\int\limits_{x_{j-1/2}}^{x_{j+1/2}}  u_{t}v dx-\int_{I_j}uv_{x}dx =
	\int\limits_{x_{j-1/2}}^{x_{j+1/2}}  uv_{x}dx - 
		\hat{u}_{j+1/2}v^{-}_{j+1/2}+\hat{u}_{j-1/2}v^{+}_{j-1/2}  \;+\\
\\
\hspace{2em}   \Bigg(C_{1} \bigg((u)^{+}_{j+1/2}-(u)^{-}_{j+1/2}\bigg)+h {C}_{2}\bigg((u_{x})^{+}_{j+1/2}-(u_{x})^{-}_{j+1/2}\bigg)\Bigg)v^{-}_{j+1/2} \; - \\
\\
\hspace{2em}  \Bigg(D_{1} \bigg((u)^{+}_{j-1/2}-(u)^{-}_{j-1/2}\bigg)+h {D}_{2}\bigg((u_{x})^{+}_{j-1/2}-(u_{x})^{-}_{j-1/2}+\bigg)\Bigg)v^{+}_{j-1/2 } \; +\\
\\
\hspace{2em}  \Bigg(h E_{1}\bigg((u)^{+}_{j+1/2}-(u)^{-}_{j+1/2}\bigg)+h^2 E_{2}\bigg((u_{x})^{+}_{j+1/2}-(u_{x})^{-}_{j+1/2}\bigg) \Bigg)(v_{x})_{j+1/2}^{-} \; -\\
\\
\hspace{2em}   \Bigg(h F_{1}\bigg((u)^{+}_{j-1/2}-(u)^{-}_{j-1/2}\bigg)+h^2 F_{2}\bigg((u_{x})^{+}_{j-1/2}-(u_{x})^{-}_{j-1/2}\bigg)\Bigg)(v_{x})_{j-1/2}^{+}
  \end{array}
  \right.
\end{equation}
\noindent where, as stated above, the following upwind flux was chosen:
\begin{equation}
\label{3_1}
\left.
\begin{array}{ll}
\hat{u}_{j+1/2}=u^{-}_{j+1/2}\;,\;\;\;\;
	\hat{u}_{j-1/2}=u^{-}_{j-1/2}
  \end{array}
  \right.
\end{equation}
\noindent The design of the above scheme is based on the minimal requirement that if $(u)_{j\pm 1/2}^{+}=(u)_{j\pm 1/2}^{-}$ and $(u_{x})_{j\pm 1/2}^{+}=(u_{x})_{j\pm 1/2}^{-}$, all penalties vanish.\\


%

\noindent We now write $u$ and $v$ in terms of their nodal values, then replacing $v(x)$ by $\varphi_{j+1/4}$ and then $\varphi_{j-1/4}$ and comparing with Eq.\eqref{eq7}, to obtain a system of eight equations for eight unknowns. The unique solution to these equations is:
%
\begin{equation}
\label{Chap4_coeff_solution_10} \nonumber
\left.
\begin{array}{llllllllll}
C_1=D_1 =-\frac{1}{2}, &&  C_2=D2=-\frac{1}{12}, && E_2= F_2=\frac{1}{72}(c_1-c_2) \\
\\
E_1=\frac{1}{36} (2c_1-6c_2+1), && &&  F_1= \frac{1}{36} (-6c_1+2c_2+1) \\
 \end{array} 
  \right.
\end{equation}
This completes the proof that the BFD scheme \eqref{eq2} can be viewed as a particular type of $p=1$, nodal-based DG scheme.





\section{Conclusions} \label{sec:Conclusions}

The BFD schemes derived for the Heat equation (see \cite{ditkowski2020error}, \cite{leblanc2020error} for further details) relied on the inherent dissipation of the diffusion operator. This dissipation, as well as the post-processing procedure, caused the damping of the high-frequency elements of the error. For the Transport problem, we had to rely on other mechanisms to manipulate the errors. 
For the correct choice of parameters, $c_{1}=c_{2}=1/2$, we could have a third-order therm, bounded in time, that can be eliminated in a post-processing stage, a fourth-order, bounded in time term, and the six-order phase error. Thus, this is a third-order or fourth-order after post-processing scheme for short to moderate times, while it is a six-order method for long times. Note that the truncation error is only of third order.

\bigskip
We also demonstrated that the BFD scheme is a particular type of nodal-based,  $p=1$,  DG scheme.

\bigskip
An immediate extension of the current work is to adapt the proposed methodology to non-periodic boundary conditions, such as Dirichlet conditions, as was done in \cite{mythesis}, for the Heat equation by creating ghost points outside the computational domain and extrapolating. Another extension would be to adapt the scheme to higher dimensions, such as 2D and 3D. The generalization to two dimensions is fairly straightforward, as the two-dimensional scheme is constructed as a tensor product of the one-dimensional scheme. This work may lead the way for constructing highly efficient DG methods.

%
%

\section{Acknowledgement}
This research was supported by Binational (US-Israel) Science Foundation grant No. 2016197.

\end{document}